\documentclass{article}

\usepackage[a4paper,left=26.5mm,right=26mm,top=26.25mm,bottom=36.75mm]{geometry}

\usepackage[numbers]{natbib}
\usepackage{amsfonts,amsmath,latexsym}
\usepackage{graphicx,xcolor,mathrsfs,mathtools}

\usepackage{hyperref}
\usepackage{authblk} 

\usepackage{amssymb}
 \usepackage{amsthm}

\newtheorem{theorem}{Theorem}
\newtheorem{corollary}[theorem]{Corollary}

\newtheorem{definition}[theorem]{Definition}

\theoremstyle{remark}

\newcommand{\R}{\mathbb{R}}
\newcommand{\N}{\mathbb{N}}

\renewcommand{\le}{\leqslant}
\renewcommand{\ge}{\geqslant}
\renewcommand{\leq}{\leqslant}

\renewcommand{\epsilon}{\varepsilon}
\renewcommand{\phi}{\varphi}

\renewcommand{\le}{\leqslant}
\renewcommand{\ge}{\geqslant}
\renewcommand{\leq}{\leqslant}

\DeclareMathAlphabet{\mathbbold}{U}{bbold}{m}{n}	

\newcommand{\xhat}{\hat{x}}
\newcommand{\xdot}{\dot{x}}
\newcommand{\D}{\mathcal{D}}
\newcommand{\UR}{\mathcal{U}}
\renewcommand{\epsilon}{\varepsilon}
\newcommand{\diff}{\mathrm{d}}
\newcommand{\fhat}{\hat{f}}

\newcommand{\length}{\ell}

\author[1]{Vincent Andrieu}
\author[1]{Lucas Brivadis}
\author[2]{Jean-Paul Gauthier}
\author[3]{Ludovic Sacchelli}
\author[1]{Ulysse Serres}

\affil[1]{Univ. Lyon, Universit\'e Claude Bernard Lyon 1, CNRS, LAGEPP UMR 5007, 43 bd du 11 novembre 1918, F-69100 Villeurbanne, France (e-mail: vincent.andrieu@univ-lyon1.fr, lucas.brivadis@univ-lyon1.fr, ulysse.serres@univ-lyon1.fr)}
\affil[2]{Universit\'e de Toulon, Aix Marseille Univ, CNRS, LIS, France (e-mail: jean-paul.gauthier@univ-tln.fr)}
\affil[3]{Department of Mathematics, Lehigh University, Bethlehem, PA, USA (e-mail: lus219@lehigh.edu)}

\title{From local to global asymptotic stabilizability for weakly contractive control systems}

\begin{document}

\maketitle

\begin{abstract}                        
A nonlinear control system is said to be weakly contractive  in the control if the flow that it generates is non-expanding (in the sense that the distance between two trajectories is a non-increasing function of time) for some fixed Riemannian metric independent of the control.
We prove in this paper that for such systems, local asymptotic stabilizability implies global asymptotic stabilizability by means of a dynamic state feedback.
We link this result and the so-called Jurdjevic and Quinn approach.
\end{abstract}


\paragraph{Acknowledgments.}   
This research was funded by the French Grant ANR ODISSE (ANR-19-CE48-0004-01).

\section{Main result}
\subsection{Statement of the result}

Consider the following nonlinear continuous-time control system:
\begin{equation}\label{eq_Syst}
\dot x=f(x,u) =f_u(x),\quad f(0,0)=0,
\end{equation}
where $x$ lives in $\R^n$ and $u$ is the control input taking values in
 an open subset $\UR$ of $\R^m$ containing zero.
We assume that $f_u\in C^1(\R^n, \R^n)$ for all $u\in\UR$,
$
\frac{\partial f}{\partial x}
\in C^0(\R^n\times\UR, \R^n)
$
and $f(x, \cdot)$ is locally Lipschitz for all $x\in\R^n$.

\begin{definition}[Static stabilizability]
System \eqref{eq_Syst} is said to be locally (resp. globally) asymptotically stabilizable by a static state feedback if there exists a locally Lipschitz mapping $\lambda:\R^n\rightarrow\UR$ such that
\begin{equation}
\dot x = f(x,\lambda(x))
\end{equation}
is locally (resp. globally) asymptotically stable
at the origin.
\end{definition}

Local asymptotic stabilizability is usually obtained by investigating first order or homogeneous approximations of the dynamical system around the origin.
Yet obtaining global stabilizability from local stabilizability is not an easy task and may fail in general.

However, there are classes of system for which we know how to bridge the gap between local and global asymptotic stabilizability.
This is obviously the case if the feedback law $\lambda$ is such that $x\mapsto f(x, \lambda(x))$ is a linear vector field.
More generally, it still holds for homogeneous systems admitting a homogeneous feedback law (see e.g. \cite{Kawski, Rosier}).
Note also that it is shown in \cite{hammouri2009two} that when the locally stabilizing state feedback fails to share the same homogeneity property than the vector field, global (or semi-global) property can still be achieved by a dynamic state feedback.

\begin{definition}[Dynamic stabilizability]
System \eqref{eq_Syst} is said to be locally (resp. globally) asymptotically stabilizable by a dynamic state feedback if there exist
$\fhat:\R^n\times\R^n\times\UR\to\R^n$
such that
$\fhat(\cdot,\cdot,u)\in C^1(\R^n\times\R^n, \R^n)$ for all $u\in\UR$,
 
$
\frac{\partial \fhat}{\partial (x, \xhat)}
\in C^0(\R^n\times\R^n\times\UR, \R^n)
$
  
and
$\fhat(x, \xhat,\cdot)$ is locally Lipschitz for all $(x, \xhat)\in\R^n\times\R^n$
and
a locally Lipschitz mapping $\lambda:~\R^n\rightarrow\UR$
such that
\begin{equation}
\dot x = f(x,\lambda(\xhat)),\quad
\dot \xhat = \fhat(x,\xhat,\lambda(\xhat))
\end{equation}
is locally (resp. globally) asymptotically stable
at the origin.
\end{definition}

In this paper, we give another class of dynamical systems which share the same property that static local asymptotic stabilizability implies dynamic global asymptotic stabilizability: namely, weakly contractive control systems.

\begin{definition}[Weakly contractive]
Let $g$ be a $C^1$ Riemannian metric on $\R^n$.
System \eqref{eq_Syst} is said to be  weakly contractive with respect to $g$ if
\begin{equation}
\forall u\in\UR, \quad L_{f_u}g \le 0,
\end{equation}
where $L_{f_u}g$ denotes the Lie derivative of the metric $g$ with respect to the vector field $f_u$.
\end{definition}

A vector field $F$ over $\R^n$ is usually said to be contractive with respect to a metric $g$ if $L_{F}g$ is negative.
Here we insist on the fact that the vector fields $f_u$ are only \emph{weakly} contractive with respect to the metric $g$, in the sense that $L_{f_u}g$ is only non-positive.

For all pair of vectors $(\phi,\psi)\in\R^n\times\R^n$, we denote by $\langle\phi,\psi\rangle$ and $|\phi|$ the canonical Euclidean inner product and induced norm over $\R^n$.
For all point $x\in\R^n$,
let $\langle \phi, \psi\rangle_{g(x)} = g(x)(\phi,\psi)$ denote the inner product between the two vectors $\phi$ and $\psi$ at the point $x$ for the metric $g$, and set $|\phi|_{g(x)} = \langle \phi, \phi\rangle_{g(x)}$.

Recall that associated to the metric $g$ we can define a distance $d_g$ between  
a pair of points
of $\R^n$ in the following way.
The length of any piecewise $C^1$ path $\gamma :[s_1,s_2]\to
\R^n$ between two arbitrary points $x_1=\gamma (s_1)$ and 
$x_2=\gamma (s_2)$ in
$\R^n$ is defined as:
\begin{equation}\label{eq_RiemanianLength}
\length(\gamma)
=\int_{s_1}^{s_2}
|\gamma'(s)|_{g(\gamma(s ))}
\diff s
\end{equation}
The distance  $d_g(x_1,x_2)$ is defined as the infimum of the length over all such paths. We denote $d_g^2$ the square of the distance function.

For all point $(x,\xhat)\in\R^n\times\R^n$, we denote (\emph{if it exists})
$\nabla_{g(\hat x)}d_g^2(x,\xhat)$
the gradient of the function $\xhat\mapsto d_g^2(x, \xhat)$ at the point $\xhat$ for the metric $g$.
 
Fix $x\in\R^n$. Then $\nabla_{g(\hat x)}d_g^2(x,\xhat)$ is well-defined if and only if, for all $\xhat\in\R^n$,
there exists a unique length-minimizing curve $\gamma$ joining $x$ to $\xhat$, \emph{i.e.} such that $\length(\gamma) = d_g(x, \xhat)$. Equivalently,
the  Riemannian exponential map at the point $\xhat$ (denoted by $\exp_{\xhat}$) is invertible\footnote{see e.g. \cite[Chap. 7, Theorem 3.1]{carmo1992riemannian} for sufficient geometric conditions.}
and we have
$$
\nabla_{g(\hat x)}d_g^2(x,\xhat)
=-2\exp_{\xhat}^{-1}(x)
$$
for all $\xhat\in\R^n$, which yields
\begin{equation}\label{eq_1}
\nabla_{g(\hat x)}d_g^2(x,\xhat) = 0
\quad\textrm{if and only if}\quad x = \xhat.
\end{equation}
Also, by definition of the Riemannian gradient, for all vectors $\phi\in\R^n$,
\begin{equation}\label{def_Grad}
\left
\langle
\nabla_{g(\hat x)}d_g^2(x,\xhat), \phi
\right
\rangle_{g(\hat x)}
=
\left\langle
\frac{\partial d_g^2}{\partial\xhat}(x,\xhat), \phi
\right\rangle.
\end{equation}

Assume that $f$ is $C^1$.
If \eqref{eq_Syst} is a  weakly contractive vector field, then for all $C^1$ control $u:\R_+\to\UR$ the time-varying vector field $f_u$ generates a non-expanding flow in the sense that,
if $x_1$ and $x_2$ satisfy $\xdot_i = f_u(x_i)$ for $i\in\{1,2\}$,
then the distance $d_g(x_1,x_2)$ between the two trajectories is a non-increasing function of time.
We give in appendix a short proof of this well-known statement to be self-contained.

The following theorem is the main result of the paper.

\begin{theorem}\label{th_main}
Let $g$ be a $C^2$ complete Riemannian metric on $\R^n$
such that $d_g^2$ is a $C^2$ function.
Assume that \eqref{eq_Syst} is  weakly contractive with respect to $g$, and $f\in C^1(\R^n\times\UR, \R^n)$.
If \eqref{eq_Syst} is locally asymptotically stabilizable by a static state feedback $\lambda\in C^1(\R^n,\UR)$,
then it is also globally asymptotically stabilizable by a dynamic state feedback
given by
\begin{equation}\label{syst_closed}
\dot x = f(x,\lambda(\xhat)),\quad
\dot \xhat = f(\xhat,\lambda(\xhat)) + k(x, \xhat)
\end{equation}
where
$$
k(x, \xhat) = -
\alpha(x, \xhat) \nabla_{g(\hat x)}d_g^2(x,\hat x)
$$
in which the function $\alpha$ has to be selected sufficiently small.
\end{theorem}
 
\subsection{Discussion on the result}
  
The idea of the proof is somehow counter-intuitive. Indeed, the feedback depends only on $\hat x$. By selecting $\alpha$ sufficiently small, we make sure that $\xhat$ remains in   the basin of attraction of the origin for the vector field associated to the state feedback. 
On the other hand, the correction terms $k$ acting on $\dot \xhat$ forces $x$ to converge to $\hat x$, which implies that $x$ goes to zero.

An interesting aspect of our approach is that no structural constraints is imposed on the local asymptotic stabilizer.
This one can be designed for qualitative purposes and can be for instance bounded or optimal as long as this one ensures a local asymptotic stability property.
This technique offers another approach to solve the global asymptotic stabilization with local optimal behavior as for instance studied in \cite{benachour2014locally} or \cite{ezal2000locally}.
The main difference with these studies being that the local optimal behavior is reproduced asymptotically in time (as $x$ converges to $\hat x$).

To construct the feedback law one needs to compute
$\nabla_{g(\hat x)}d_g^2(x,\hat x)$ which may be difficult to obtain analytically in general
(except in some simple cases, \emph{e.g.}, if the metric is constant).
Some ways of constructing similar correction terms may be obtained following observer designs based on Riemannian approaches as in \cite{Aghannan,sanfelice2011convergence}.
In particular in \cite[Lemma 3.6]{sanfelice2011convergence}, the authors introduced a ``distance-like'' function $\delta$, that is of crucial importance in the construction of the correction term.

\subsection{Proof}
Let $\lambda$ be a $C^1$ locally asymptotically stabilizing feedback law.
Let $\D$ be the basin of attraction of the origin for the vector field $x\mapsto f(x, \lambda(x))$, which is a non-empty open subset of $\R^n$.
 
According to the converse Lyapunov theorem \cite{teel2000smooth} (based on the previous works of \cite{kurzweiloriginal, kurzweil, massera})  , there exists a proper function $V\in C^\infty(\D,\R_+)$ such that $V(0)=0$ and
 
\begin{equation}\label{eq_Lyap}
\frac{\partial V}{\partial x}(x)f(x,\lambda(x)) \le -V(x),\quad \forall x\in\D\ .
\end{equation}
For all $r>0$, set $D(r) = \{x\in\R^n\mid V(x) \le r\} $ which is a compact subset of $\D$.
Let $\alpha:\R^n\times\D\to\R_+$ be the positive and locally Lipschitz function given by
\begin{equation}
\alpha(x, \xhat) = 
\frac{-\max\{V(\hat x),1\}}{2\left(1 + \left|\frac{\partial V}{\partial x}(\hat x)\right|\right)
\left(1+
\left|\nabla_{g(\hat x)}d_g^2(x,\hat x)\right|\right)}
.
\end{equation}
It yields
\begin{equation}\label{E:norm_k}
|k(x,\xhat)| \le \frac{\max\{V(\hat x),1\}}{2\left(1 + \left|\frac{\partial V}{\partial x}(\hat x)\right|\right)},\quad \forall(x, \xhat)\in\R^n\times\D.
\end{equation}

We prove Theorem~\ref{th_main} in three steps.

\noindent\textbf{Step 1 : 
the $\hat x$-component of semi-trajectories of \eqref{syst_closed} remain in a compact subset of $\D$.}
For all $(x, \hat x)\in\R^n\times\D$,
it follows from \eqref{eq_Lyap} and \eqref{E:norm_k} that
$$
\begin{aligned}
\frac{\partial V}{\partial x}(\xhat)[f(\hat x,\lambda(\hat x)) + k(x,\hat x)]
&\le  -V(\hat x) + \left|\frac{\partial V}{\partial x}(\xhat)\right| \frac{\max\{V(\hat x),1\}}{2\left(1 + \left|\frac{\partial V}{\partial x}(\xhat)\right|\right)}
\\
&\le  -V(\hat x) + \frac{1}{2}\max\{V(\hat x),1\}.
\end{aligned}
$$
Hence, if $\xhat\in \D\setminus D(1)$,
\begin{equation}
    \frac{\partial V}{\partial x}(\xhat)(f(\hat x,\lambda(\hat x)) + k(x,\hat x))\le -\frac{1}{2} V(\hat x)\ .
\end{equation}
For all initial conditions $(x_0,\hat x_0)\in\R^n\times\D$, the 
solution $(x, \hat x)$ of the closed-loop system \eqref{syst_closed}  satisfies
\begin{equation*}
V(\hat x(t))\le \max\{ V(\hat x_0), 1\},
\end{equation*}
for all $t\ge0$, in the time domain of existence of the solution.
In other words, $\hat x(t)\in D(1)\cup D(V(\hat x_0))$ which is a compact subset of $\D$.

\noindent\textbf{Step 2 : the distance between $\hat x$ and $x$ is non-increasing and 
has limit zero.}
System \eqref{syst_closed} can be rewritten as
\begin{equation}
    \begin{bmatrix}
    \dot x\\
    \dot \xhat
    \end{bmatrix}
    = F(x,\hat x) + K(x,\hat x)
\end{equation}
by setting
$
F(x,\hat x) = \begin{bmatrix}
    f(x,\lambda(\xhat))\\
    f(\xhat,\lambda(\xhat))
    \end{bmatrix}$
and
$K(x,\hat x)= \begin{bmatrix}
    0\\
    -\alpha(x,\xhat)\nabla_{g(\hat x)}d_g^2(x, \xhat)
    \end{bmatrix}
$.

Since \eqref{eq_Syst} is  weakly contractive with respect to $g$, the result proved in appendix applied to the control $u = \lambda(\xhat)$ shows that
$$
L_Fd_g^2(x,\xhat)\le 0.
$$
Thus, by \eqref{def_Grad},
\begin{equation}\label{eq_dist}
L_{F+K}d_g^2(x,\hat x) \le -\alpha(x,\xhat)  \left|\nabla_{g(\hat x)}d_g^2(x, \xhat)\right|_{g(\hat x)}^2.
\end{equation}
Hence, for all $(x_0,\hat x_0)\in\R^n\times\D$, $t\mapsto d_g(x(t),\hat x(t))$ is non-increasing and for all $t\ge 0$ on the time domain of existence of the solution we have
$$
(x(t),\hat x(t)) \in \Gamma(x_0,\hat x_0) ,
$$
where
$$
\Gamma(x_0,\hat x_0) = \Big\{(\xi,\hat \xi)\in\R^n\times\D\mid \hat \xi\in D(1)\cup D(V(\hat x_0)), 
d_g(\xi,\hat \xi)\le d_g(x_0,\hat x_0)\Big\} .
$$
Moreover,
$g$ is a complete metric. Then, according to the Hopf-Rinow theorem,
$\Gamma(x_0,\hat x_0)$ is compact.
Hence, solutions of \eqref{syst_closed} are complete in positive time.

Given $(x_0,\hat x_0) \in \R^n\times\D$, let $\kappa:\R_+\to\R_+$ be the function defined by
$$
\kappa(s) = \min_{(\xi,\hat \xi)\in\Gamma(x_0,\hat x_0)\mid d_g(\xi,\hat \xi)=s}\alpha(\xi,\hat \xi)\left|\nabla_{g(\hat \xi)}d_g^2(\xi, \hat \xi)\right|_{g(\hat \xi)}^2.
$$
Note that
if $x_0 \neq \hat x_0$, then,
for all $s>0$, $\kappa(s)>0$ since $\alpha$ takes positive values and \eqref{eq_1} holds.  
Hence, \eqref{eq_dist} leads to
\begin{equation}\label{eq_dist2}
\frac{\diff}{\diff t} d_g^2(x(t),\hat x(t)) \le -\kappa(d_g^2(x(t),\hat x(t)) )\ , \ \forall t\ge 0.
\end{equation}
Thus
$\lim_{t\to+\infty}d_g(x(t),\hat x(t))=0$.

\noindent\textbf{Step 3 : attractivity and local asymptotic stability of the origin.}
Given $(x_0,\hat x_0)$ in $\R^n\times\D$, let $\mu:\R_+\to\R_+$ be the function defined by
$$
\mu(s) = \max_{(\xi, \hat \xi)\in\Gamma(x_0,\hat x_0)\mid d_g(\xi,\hat \xi)\le s}\left|\frac{\partial V}{\partial x}(\hat \xi)k(\xi,\hat \xi)\right|.
$$
Then $\mu$ is non-decreasing, continuous and $\mu(0)=0$.
Moreover,
the solution $(x,\hat x)$ of \eqref{syst_closed} initialized at $(x_0,\hat x_0)\in\R^n\times\D$ satisfies
\begin{equation}\label{eq_2}
 \frac{\diff }{\diff t}V(\xhat(t)) \le -V(\hat x(t)) + \mu (d_g(x(t),\hat x(t)).
\end{equation}
From this inequality and Step 2 we conclude that $\lim_{t\rightarrow +\infty}(x(t),\hat x(t))=(0,0)$.

Inequalities \eqref{eq_dist} and \eqref{eq_2} being true for all solutions starting in $\Gamma(x_0,\hat x_0)$, this implies also stability of $(0, 0)$.

\section{Link with Jurdjevic and Quinn approach}

\subsection{Jurdjevic and Quinn result}

The next result follows from the work of Jurdjevic and Quinn in \cite{jq}. The version that we state here is a direct corollary of \cite[Theorem II.1]{mazenc}

\begin{theorem}[Jurdjevic and Quinn approach]\label{th_jq}
Consider the control system
\begin{equation}\label{syst_jq}
    \dot x = a(x) + b(x,u)u,
\end{equation}
with $a$ and $b$ two $C^1$ functions.
Assume that there exists a $C^1$ positive definite proper function
$V:\R^n\mapsto\R_+$ such that
$$
L_aV\leqslant 0.
$$
If the only solution of the system
\begin{equation}
\dot x = a(x),\quad
L_{b(\cdot, 0)}V(x)=0,\quad
L_{a}V(x)=0
\end{equation}
is $x\equiv0$,
then \eqref{syst_jq} is globally asymptotically stabilizable by a static state feedback.
\end{theorem}

In the context of weakly contractive control systems, the Jurdjevic and Quinn approach leads to the following corollary.

\begin{corollary}\label{cor_jq}
Let $g$ be a complete Riemannian metric on $\R^n$.
Assume that \eqref{eq_Syst} is  weakly contractive with respect to $g$ and that $f\in C^2(\R^n\times\UR,\R^n)$.
If the only solution of the system
\begin{equation*}
\dot x = f(x,0),\
\left(L_{b(\cdot, 0)}d_g^2(\cdot,0)\right)(x)=0,\
L_{f_0}d_g^2(x,0)=0
\end{equation*}
where $b(\cdot, 0) = \frac{\partial f}{\partial u}(\cdot,0)$
is $x\equiv0$,
then \eqref{eq_Syst} is globally asymptotically stabilizable by a static state feedback.
\end{corollary}

To prove this corollary, it is sufficient to apply Theorem~\ref{th_jq} with $V:x\mapsto d_g(x,0)^2$, $a:x\mapsto f(x, 0)$ and $b:(x, u)\mapsto \int_0^1\frac{\partial f}{\partial u}(x,su)ds$.

\subsection{Link with our result}

Note that the Jurdjevic-Quinn approach guarantees the existence of a \emph{static} state feedback, contrarily to our main Theorem \ref{th_main} which build a \emph{dynamic} state feedback.
However, the feedback obtained by their approach is implicit, while our dynamic state feedback is explicitly given by \eqref{syst_closed}.

Moreover, our feedback law differs strongly with the one given in Jurdjevic-Quinn approach. Indeed, in their approach the feedback is designed small enough to make sure that it acts in a good direction related to the Lyapunov function. 
In our framework, this is no more a \textit{small feedback approach} but more a \textit{small correction term for an observer approach}.

Let us consider the particular case in which
\begin{equation}
\label{syst_example}
f(x,u) = Ax + b(x)u.
\end{equation}
where $A\in\R^{n\times n}$ and $b\in C^1(\R^n, \R^{n})$.
Then \eqref{syst_example} is  weakly contractive with respect to some constant metric $g$ if and only if
$L_Ag \leq 0$ and $L_bg = 0$\footnote{It is easy to check that it is the case if and only if $b(x) = b(0) + Jx$ with $L_Jg = 0$.}.
Moreover, the pair $(A, b(0))$ is controllable if and only if \eqref{syst_example} is locally asymptotically stabilizable by a static feedback.
Then, if all these hypotheses hold, a dynamic globally stabilizing state feedback is given by Theorem~\ref{th_main}.

We can also show under the same hypotheses that the Jurdjevic and Quinn approach can be applied.
Indeed, the system in Corollary~\ref{cor_jq} is equivalent to
\begin{equation}
\dot x = Ax,\quad
\left(L_{b}d_g^2(\cdot,0)\right)(x)=0,\quad
L_{A}d_g^2(x, 0)=0
\end{equation}
which implies that $x\equiv0$
when the pair $(A,  b(0))$ is controllable.
Then, according to Corollary~\ref{cor_jq}, \eqref{syst_example} is globally asymptotically stabilizable by a static state feedback.
{
However, it is not clear in general that both contexts are equivalent, and finding an example fitting in the framework of Theorem \ref{th_main} but in which the Jurdjevic and Quinn approach of Corollary~\ref{cor_jq} remains an open question.
}

\section{Appendix on weakly contractive vector fields}

For all $u:\R_+\to\UR$ and all $x\in\R^n$, denote by $t\mapsto X_u(x, t)$ the solution of \eqref{eq_Syst} with initial condition $x$.
Let $u:\R_+\to\UR$ be such that $X_u$ is well-defined and $C^2$ on $\R^n\times\R_+$.

Let $(x_1, x_2)\in\R^n\times\R^n$
and $\gamma:[s_1,s_2]\to
\R^n$ be a $C^2$ path between the points $x_1=\gamma (s_1)$ and 
$x_2=\gamma (s_2)$.
For all $(s, t)\in[s_1, s_2]\times\R_+$,
set $\Gamma(s, t) = X_u(\gamma(s), t)$
and $\rho(s, t) = \left|\frac{\partial\Gamma}{\partial s}(s, t)\right|^2_{g(\Gamma(s, t))}$.
Then $\rho$ is $C^1$ and
\begin{align*}
    \frac{\partial \rho}{\partial t}(s, t)
    =L_{f_u}g(\Gamma(s, t))\left(\frac{\partial\Gamma}{\partial s}(s, t), \frac{\partial\Gamma}{\partial s}(s, t)\right)
    \leq 0,
\end{align*}
which yields
\begin{align*}
\frac{\diff \length(\Gamma(\cdot, t))}{\diff t}
&=\frac{\diff}{\diff t} \int_{s_1}^{s_2}
\sqrt{\rho(s, t)}
\diff s\\
&= \int_{s_1}^{s_2}
\frac{1}{2\sqrt{\rho(s, t)}}
\frac{\partial \rho}{\partial t}(s, t)
\diff s\\
&\leq 0.
\end{align*}
Hence
$
d_g(X_u(x_1, t), X_u(x_2, t))
\leq \length(\Gamma(\cdot, t))
\leq \length(\gamma).
$
Choosing a sequence of paths $(\gamma_n)_{n\in\N}$ such that $\length(\gamma_n)\to d_g(x_1, x_2)$ and passing to the limit we get
\begin{align*}
d_g(X_u(x_1, t), X_u(x_2, t))
\leq d_g(x_1, x_2).
\end{align*}
Since this inequality is true for any control input $u$, $t\mapsto d_g(X_u(x_1, t), X_u(x_2, t))$ is non-increasing for all control $u$ and all points $x_1$, $x_2$.

\bibliographystyle{abbrv}
\bibliography{references}

\end{document}